\documentclass[12pt]{article}

\usepackage{amssymb,amsmath,amsthm}

\newcommand{\vp}{\varepsilon}
\newcommand{\cl}[1]{{\mathcal{#1}}}
\newcommand{\bb}[1]{{\mathbb{#1}}}
\newcommand{\tr}{\text{tr}}

\theoremstyle{plain}
\newtheorem{pro}{Proposition}[section]
\newtheorem{lem}[pro]{Lemma}
\newtheorem{cor}[pro]{Corollary}
\newtheorem{thm}[pro]{Theorem}

\theoremstyle{definition}
\newtheorem{defn}[pro]{Definition}

\theoremstyle{remark}
\newtheorem{rem}[pro]{Remark}

\begin{document}
\title{STRONGLY SINGULAR MASAS\\ IN TYPE ${\mathrm {II}}_1$ FACTORS}

\author{Allan M.\ Sinclair\\
Department of Mathematics\\
University of Edinburgh\\
Edinburgh EH9 3JZ\\
SCOTLAND\\
{\tt e-mail:\ allan@maths.ed.ac.uk}
\and
Roger R. Smith\footnote{Partially supported by a grant from the National Science 
Foundation.}\\
Department of Mathematics\\
Texas A\&M University\\
College Station, TX \ 77843\\ 
U.S.A.\\
{\tt e-mail:\ rsmith@math.tamu.edu}}

\date{}

\maketitle{}

\bigskip
\begin{abstract}
In this paper we introduce and study strongly singular maximal abelian
self--adjoint subalgebras of type $II_1$ factors. We show that certain elements
of free groups and of non--elementary hyperbolic groups generate such masas,
and these also give new examples of masas for which Popa's invariant
$\delta(\cdot)$ is $1$. We also explore the connection between Popa's invariant
and strong singularity.  \end{abstract} 
\newpage 
\setcounter{section}{0}
\section{Introduction}\label{sec1} \indent

The study of maximal abelian self-adjoint subalgebras (masas) in a von
Neumann algebra ${\mathcal{M}}$ has a long and rich history,
\cite{B,Dix,Dyk,Ge,GP,JP,KS,Po1,Po2,Pu1,Pu2}. Various types of masas have been
identified and investigated, often categorized by their groups of normalizing
unitaries. At one end of the spectrum are the regular or Cartan masas; these
have sufficiently many normalizing unitaries to generate ${\mathcal{M}}$. At
the other end are the singular masas; the only normalizing unitaries of such a
masa ${\mathcal{A}}$ are the unitaries of  ${\mathcal{A}}$, \cite{Dix}. In this paper we
introduce a new class of singular masas, which we call strongly singular.
These are defined by an inequality relating the distance between
${\mathcal{A}}$ and a unitary conjugate $u{\mathcal{A}}u^*$ to the distance of
$u$ to ${\mathcal{A}}$, also allowing us to introduce a new invariant $\alpha({\cl A})$
for masas, taking values in $[0,1]$
(definitions are contained in the second section). Part of our original 
motivation was the observation that a
reverse inequality between these quantities is always valid (Proposition
\ref{pro2.1}). 

From the point of view of the inner automorphism group, a Cartan masa ${\cl A}$
is flexible in that any two projections in ${\cl A}$ with equal trace can be switched 
by an inner automorphism of ${\cl M}$ which leaves ${\cl A}$ invariant, \cite{Po1}.
For a singular masa, any inner automorphism of ${\cl M}$ which leaves ${\cl A}$ invariant 
has trivial action on the masa, but this takes little account of the other unitaries in
${\cl M}$. Popa's invariant $\delta({\cl A})$ of a masa ${\cl A}$, 
\cite{Po2}, is a measure 
of this rigidity in terms of partial isometries in ${\cl M}$ whose 
initial and final projections are orthogonal in ${\cl A}$. Strong singularity is intended to 
develop a rigidity condition on ${\cl A}$ which reflects the perturbations of the masa by 
the inner automorphisms of ${\cl M}$. This condition is compatible with the 
isometric action of the inner automorphism group of ${\cl M}$ on the natural metric space 
of the masas in ${\cl M}$.

We now describe the contents of the paper. The second section contains 
definitions and some preliminary  results, while the
third section presents some examples of strongly singular masas. The
main results here are that the masas arising from the generators of free
groups (Corollary \ref{cor3.4}) and, more generally, from prime elements of
Gromov's non--elementary I.C.C. hyperbolic groups, \cite{Gr}, (Theorem
\ref{thm3.6}) are strongly singular. The {\it {Infinite Conjugacy Class}}
condition is included to ensure that the resulting von Neumann algebras are
factors, \cite[p.126]{Ha}. This condition for a non-elementary hyperbolic group is
equivalent to the group being torsion free.
The techniques also show that these masas satisfy
$\delta({\mathcal{A}})=1$.  
The singular masa in the hyperfinite
type $II_1$ factor constructed by Tauer, \cite{Ta}, was the only previous
example where $\delta(\cdot)$ could be exactly determined, \cite{Po2}, although
Popa had shown that every type $II_1$ factor contains a masa with
$\delta(\cdot) \geq 10^{-4}$, \cite{Po2}. By tightening the argument given by
Popa, we can improve this bound to $1/58$, but it seems difficult to obtain
any estimate close to 1. We also give an example of a strongly singular masa
in the hyperfinite type $II_1$ factor whose Popa invariant is 1 (Corollary
\ref{cor3.8}).

   Our examples of strongly
singular masas have conditional expectations which satisfy a multiplicative
condition, which we use to define an asymptotic homomorphism in the fourth
section. The main result (Theorem \ref{thm6.7}) is that masas whose
conditional expectations are asymptotic homomorphisms all satisfy
$\delta({\mathcal{A}})=1$.
In the last section of the paper, we relate strong singularity to Popa's
invariant, and we prove that every masa satisfies
$\delta({\mathcal{A}}) \geq \alpha({\cl A})/\sqrt{5}$, (Theorem \ref{thm7.1}).
In particular,  $\delta({\mathcal{A}}) \geq 1/\sqrt{5}$ for strongly singular
masas. Intuitively, it would seem reasonable that these two invariants should
be equal, or at  least mutually dominating. However, we have been unable to
obtain a reverse inequality of the form $\delta({\mathcal{A}}) \leq c \cdot
\alpha({\cl A})$ for some constant $c > 0$. 

We conclude by mentioning that a significant part of our work has been
motivated by the papers of Sorin Popa on masas, \cite{Po1,Po2,Po3,Po4}, and
particularly by the results on orthogonality in \cite{Po3}. We also
thank Pierre de la Harpe for pointing out an error in an earlier
version of the paper. \newpage      

\section{Preliminaries}\label{sec2}
\setcounter{equation}{0}
\indent

In this section we present some definitions and notation which we will use
subsequently. In order to motivate the definition of an $\alpha$-strongly
singular masa, we will first prove an easy inequality concerning unitary
conjugates of subalgebras.

Throughout we denote the operator norm on a type $II_1$ factor $\cl M$ by
$\|\cdot\|$, while $\|\cdot\|_2$ denotes the norm $(\tr~(x^*x))^{1/2}$ induced
by the unique normalized trace. We write $L^2(\cl M,\tr)$ for the Hilbert space
completion of $\cl M$ in $\|\cdot\|_2$. A linear map $\phi\colon \ \cl M\to \cl M$ may be
viewed as having range in $L^2(\cl M, \tr)$. If it is then bounded, we denote
its norm by $\|\phi\|_{\infty,2}$. If it is also bounded as a map on $L^2(\cl
M, \tr)$, we write $\|\phi\|_2$ for this norm. We reserve $\|\phi\|$ for the
norm when $\cl M$ has the operator norm for both range and domain. For each
von~Neumann subalgebra $\cl N$, there is a unique trace preserving conditional
expectation $\bb E_{\cl N} \colon \ \cl M\to \cl N$, and it is contractive for
each of the norms $\|\cdot\|$, $\|\cdot\|_2$ and $\|\cdot\|_{\infty,2}$.
Moreover, if $\phi,\psi\colon \ \cl M\to \cl M$ are linear maps, then the
inequalities \begin{equation}\label{eq2.1}
\|\phi\psi\|_{\infty,2} \le \|\phi\|_{\infty,2} \|\psi\|,\quad \|\phi\|_2
\|\psi\|_{\infty,2}
\end{equation}
are immediate from the definitions. We note, for future reference, one
important property of conditional expectations: \ $\bb E_{\cl N}$ is an $\cl
N$-bimodule map, \cite{Tom}.

\begin{pro}\label{pro2.1}
Let $\cl M$ be a type $II_1$ factor and let $\cl A$ be a von Neumann
subalgebra. For any unitary $u\in \cl M$,
\begin{equation}\label{eq2.2}
\|\bb E_{\cl A} - \bb E_{u\cl A u^*}\|_{\infty,2} \le 4\|u-\bb E_{\cl
A}(u)\|_2.
\end{equation}
\end{pro}

\begin{proof}
For any $x\in \cl M$,
\begin{equation}\label{eq2.3}
\bb E_{u\cl A u^*} (x) = u\bb E_{\cl A}(u^*xu)u^*,
\end{equation}
and so
\begin{align}
\|(\bb E_{\cl A} - \bb E_{u\cl Au^*})(x)\|_2 &= \|\bb E_{\cl A}(x) - u\bb
E_{\cl A}(u^*xu)u^*\|_2\nonumber\\
\label{eq2.4}
&= \|\bb E_{\cl A}(x)u - u\bb E_{\cl A}(u^*xu)\|_2.
\end{align}
To estimate $\|\bb E_{\cl A}-\bb E_{u\cl Au^*}\|_{\infty,2}$, we may replace
$x$, $\|x\| \le 1$, in (\ref{eq2.4}) by unitaries $w$, and we may further
assume that $w$ has the form $uv$ for some unitary $v\in \cl M$. Thus it
suffices to estimate
\begin{equation}\label{eq2.5}
\|\bb E_{\cl A}(uv)u - u\bb E_{\cl A}(vu)\|_2
\end{equation}
as $v$ ranges over the unitary group of $\cl M$.

Write $a = \bb E_{\cl A}(u) \in \cl A$, $b = (I-\bb E_{\cl A})(u)$. Then $u =
a+b$, and $\|a\|^2_2 + \|b\|^2_2 = 1$. Thus (\ref{eq2.5}) becomes
\begin{align}
&\|\bb E_{\cl A}(av+bv)u - u\bb E_{\cl A}(va+vb)\|_2\nonumber\\
&\quad \le \|\bb E_{\cl A}(bv)u\|_2 + \|u\bb E_{\cl A}(vb)\|_2 + \|\bb E_{\cl
A}(av)(a+b)-(a+b)\bb E_{\cl A}(va)\|_2\nonumber\\
&\quad\le 2\|b\|_2 + \|\bb E_{\cl A}(av)b\|_2 + \|b\bb E_{\cl A}(va)\|_2 +
\|a\bb E_{\cl A}(v) a-a\bb E_{\cl A}(v)a\|_2\nonumber\\
&\quad \le 4\|b\|_2\nonumber\\
\label{eq2.6}
&\quad = 4\|u-\bb E_{\cl A}(u)\|_2.
\end{align}
The result follows by taking the supremum over all unitaries $v$ in
(\ref{eq2.6}).
\end{proof}

It is natural to ask whether a reverse inequality of the form
\begin{equation}\label{eq2.7}
\|\bb E_{u\cl Au^*} - \bb E_{\cl A}\|_{\infty,2} \ge \alpha\|u-\bb E_{\cl
A}(u)\|_2
\end{equation}
can hold for some $\alpha>0$ and for all unitaries $u\in \cl M$. Of course
some restrictions on this question must be made, because (\ref{eq2.7}) forces
any normalizing unitary of $\cl A$ to lie in $\cl A$. This rules out abelian
algebras which are not maximal, regular and semi-regular masas, and any
algebra $\cl A$ for which $\cl A'\not\subseteq \cl A$. Thus any masa which
satisfies (\ref{eq2.7}) is automatically singular. We will show subsequently
that many singular masas satisfy such an inequality, and this suggests the
following terminology.

\begin{defn}\label{defn2.2}
A masa $\cl A$ in a type $II_1$ factor $\cl M$ is said to be $\alpha$-strongly
singular if (\ref{eq2.7}) holds. When $\alpha=1$, we say that $\cl A$ is
strongly singular. We let $\alpha(\cl A)$ denote the supremum of all numbers
$\alpha$ for which (\ref{eq2.7}) is valid. We note that strong
singularity and $\alpha(\cdot)$ can be defined in this way for
any von~Neumann subalgebra. $\hfill\square$  
\end{defn}

It is clear, from Proposition \ref{pro2.1}, that $\alpha(\cl A)$ takes its
value in [0,4]. Since we will construct examples where $\alpha(\cl A)= 1$,
the following result gives the optimal upper estimate on $\alpha(\cdot)$.

\begin{pro}\label{pro2.3} If $\cl A$ is a masa in a type $II_1$ factor $\cl M$,
then $\alpha(\cl A) \leq 1$. \end{pro}

\begin{proof} Let $u \in \cl M$ be any unitary, and regard $\bb E_{\cl A}$ and 
$\bb E_{u\cl Au^*}$ as projections in $B(L^2({\cl M}),\tr)$. Both are positive operators, so
the inequality
\begin{equation}\label{eq2.7a}
\|\bb E_{u\cl Au^*} - \bb E_{\cl A}\|_2 \leq 
{\mathrm {max}}\,\{\|\bb E_{u\cl Au^*}\|_2,\,\|\bb E_{\cl A}\|_2\}=1
\end{equation}
follows by applying states to this difference. 
Then $\|\bb E_{u\cl Au^*} - \bb E_{\cl A}\|_{\infty,2}\leq 1$, since 
$\|\cdot\|_2 \geq \|\cdot\|_{\infty,2}$ by taking $\psi = I$ in (\ref{eq2.1}).

Now choose a projection $p \in \cl A$, $\tr(p)=1/2$, and choose a partial isometry
$v \in \cl M$ such that
\begin{equation}\label{eq2.7b}
vv^*=p,\ \ v^*v=p^{\perp}.
\end{equation}
The element $u=v+v^*$ is a unitary in $\cl M$ satisfying $upu=p^{\perp}$ or, 
equivalently, $pu=up^{\perp}$. Then $p\bb E_{\cl A}(u)=\bb E_{\cl
A}(u)p^{\perp}$, which forces $\bb E_{\cl A}(u)=0$, since these operators
commute. Thus \begin{equation}\label{eq2.7c}
\|u-\bb E_{\cl A}(u)\|_2=\|u\|_2=1.
\end{equation}
This choice of unitary shows that the inequality in (\ref{eq2.7}) fails for each $\alpha > 1$, and it follows that 
$\alpha(\cl A) \leq 1$.
\end{proof}

Let $\cl A$ be a masa in a type $II_1$ factor $\cl M$, and let $v$ be a
non-zero partial isometry in $\cl M$ such that $p=vv^*$ and $q=v^*v$ are
orthogonal projections in $\cl A$. Define $\delta(v\cl Av^*,\cl A)$ by 
\begin{equation}\label{eq2.8}
\delta(v\cl Av^*,\cl A) = \sup\{\|x-\bb E_{\cl A}(x)\|_2\colon \ x\in v\cl A
v^*,\ \|x\|\le 1\}.
\end{equation}
Then $\delta(\cl A)$ is the largest number $\lambda$ for which the inequality
\begin{equation}\label{eq2.9}
\delta(v\cl A v^*,\cl A) \ge \lambda\|v^*v\|_2
\end{equation}
holds for all such partial isometries (see \cite{Po2}). Since any element $x\in
v\cl Av^*$ satisfies $x=pxp$, it is clear that
\begin{align}
\delta(v\cl Av^*,\cl A)^2 &\le \sup\{\|pxp\|^2_2\colon \ x\in \cl M,
\ \|x\|\le 1\}\nonumber\\
\label{eq2.10}
&\le \tr(p) = \tr(q) = \|q\|^2_2 = \|v^*v\|^2_2.
\end{align}
It follows from (\ref{eq2.10}) that $\delta(\cl A)\le 1$. This was stated in
\cite{Po2}, but we have included a proof for the reader's convenience. \newpage

\section{Strong singularity in discrete group factors}\label{sec3}

\setcounter{equation}{0}

\indent
In this section we present some examples of strongly singular masas and
subfactors arising from discrete groups. When $\Gamma$ is a discrete I.C.C.
group (each element other than the identity has an infinite conjugacy class)
the resulting von~Neumann algebra $VN(\Gamma)$, represented on
$\ell^2(\Gamma)$, is a type $II_1$ factor. Each element of the group is a
unitary in $VN(\Gamma)$ and thus generates an abelian von~Neumann subalgebra.
Since the principal examples of type $II_1$ factors arise from discrete
groups, our examples of strongly singular masas will be generated by elements
of groups. The first two lemmas give key technical results which will be needed
for our main theorems. The common hypotheses for the first three results are
taken from \cite{Po3}. Note that the I.C.C. hypothesis is inessential for the
proofs, and is only included to ensure that the associated von~Neumann algebras
are factors.

\begin{lem}\label{lem3.1}
Let $G$ be an infinite subgroup of a countable discrete I.C.C. group $\Gamma$ 
with the property that $xGx^{-1}\cap G=\{e\}$ for all $x\in \Gamma \backslash
G$, let ${\cl M}=VN(\Gamma)$, and let ${\cl N}=VN(G)$. For each set of elements
$u_i \in {\bb C}\Gamma$, $1 \leq i \leq n$, the equations 
\begin{equation}\label{eq3.1}
{\bb E}_{\cl N}(u_sgu_t)={\bb E}_{\cl N}(u_s)g{\bb E}_{\cl N}(u_t),\ \ 1 \leq
s,t \leq n,
\end{equation}
are satisfied by all but a finite number of $g \in G$.
\end{lem}

\begin{proof}
Each $u_i$ is a finite linear combination of group elements, so it suffices to
prove, for  fixed $h,k \in \Gamma$, that the equation
\begin{equation}\label{eq3.2}
{\bb E}_{\cl N}(hgk)={\bb E}_{\cl N}(h)g{\bb E}_{\cl N}(k)
\end{equation}
is satisfied by all but a finite number of $g \in G$. The modular properties
of  ${\bb E}_{\cl N}$ show that (\ref{eq3.2}) always holds when either $h$ or
$k$ is in $G$, so we may assume that both elements are not. In this case the
right hand side of (\ref{eq3.2}) is 0, and so we only need establish that $hgk
\in G$ for only finitely many $g \in G$. If $g_1$ and $g_2$ are two such
elements, then $hg_1g_2^{-1}h^{-1} \in G$. The hypotheses then imply that 
$g_1=g_2$, showing that (\ref{eq3.2}) fails for at most one $g \in G$.
\end{proof}

\begin{lem}\label{lem3.2}
Let $G$ be an infinite subgroup of a countable discrete I.C.C. group $\Gamma$ 
with the property that $xGx^{-1}\cap G=\{e\}$ for all $x\in \Gamma \backslash
G$, let ${\cl M}=VN(\Gamma)$, and let ${\cl N}=VN(G)$. If $u,v \in {\cl M}$, then 
\begin{align}
\|(I-\bb E_{\cl N})(u\bb E_{\cl N}(\cdot)v)\|^2_{\infty,2} &\ge 
\tr[\bb E_{\cl N}(\bb E_{{\cl N}'\cap {\cl M}}(u^*u)) 
(\bb E_{\cl N}(vv^*) -  \bb E_{\cl N}(v) \bb E_{\cl N}(v)^*)]\nonumber\\
\label{eq3.3}&\quad -\|(I-\bb E_{{\cl N}'\cap {\cl M}})(u^*u)\|_2
(\|vv^*\|_2+\|\bb E_{{\cl N}}(v)\bb E_{{\cl
N}}(v)^*\|_2). \end{align}
If $u$ is a unitary, then
\begin{equation}\label{eq3.4}
\|(I-\bb E_{\cl N})(u\bb E_{\cl N}(\cdot)v)\|^2_{\infty,2} \ge
\tr[\bb E_{\cl N}(vv^*) -  \bb E_{\cl N}(v) \bb E_{\cl N}(v)^*].
\end{equation}
\end{lem}

\begin{proof}
If we can prove (\ref{eq3.3}) for $u,v \in {\bb C}\Gamma$, then the 
$\|\cdot \|_2$--norm continuity of conditional expectations and the Kaplansky
density
theorem will show that it holds generally. Thus
we assume that $u,v \in {\bb C}\Gamma$. Then, by Lemma \ref{lem3.1},
we may choose $g \in G$ so that \begin{equation}
\label{eq3.5}
{\bb E}_{\cl N}(ugv)={\bb E}_{\cl N}(u)g{\bb E}_{\cl N}
(v),\ \ {\bb E}_{\cl N}(u^*ugvv^*)={\bb E}_{\cl N}(u^*u)g
{\bb E}_{\cl N}(vv^*).
\end{equation}
For this choice of $g$,
\begin{align}
\|(I-\bb E_{\cl N})(u\bb E_{\cl N}(\cdot)v)\|^2_{\infty,2} &\ge \|(I-\bb
E_{\cl N})(ugv)\|^2_2\nonumber\\
&= \|ugv\|^2_2 - \|\bb E_{\cl N}(ugv)\|^2_2\nonumber\\
&= \tr(u^*ugvv^*g^{-1}) - \|\bb E_{\cl N}(ugv)\|^2_2\nonumber\\
&= \tr(\bb E_{\cl N}(u^*ugvv^*g^{-1}) - \|\bb E_{\cl N}(ugv)\|^2_2\nonumber\\
\label{eq3.6}&= \tr(\bb E_{\cl N}(u^*u)g\bb E_{\cl N}(vv^*)g^{-1}) - \|\bb
E_{\cl N}(ugv)\|^2_2, \end{align}
where the last equality follows from (\ref{eq3.5}), since $g^{-1} \in \cl N$.
Now write \begin{equation}\label{eq3.7}
a=\bb E_{{\cl N}'\cap {\cl M}}(u^*u),\ \ b=(I-\bb E_{{\cl
N}'\cap {\cl M}})(u^*u).
\end{equation}
Since $ga=ag$, we may apply $\bb E_{{\cl N}}$ to conclude that $g$ and $\bb
E_{{\cl N}}(a)$ commute. Thus
 \begin{align}
\tr (\bb E_{{\cl N}}(u^*u)g\bb E_{{\cl N}}(vv^*)g^{-1})
&=
\tr (\bb E_{{\cl N}}(a+b)g\bb E_{{\cl N}}(vv^*)g^{-1})\nonumber\\
&=
\tr (g\bb E_{{\cl N}}(a)\bb E_{{\cl N}}(vv^*)g^{-1})+
\tr (\bb E_{{\cl N}}(b)g\bb E_{{\cl N}}(vv^*)g^{-1})\nonumber\\
&=
\tr (\bb E_{{\cl N}}(a)\bb E_{{\cl N}}(vv^*))+
\tr (\bb E_{{\cl N}}(b)g\bb E_{{\cl N}}(vv^*)g^{-1})\nonumber\\
\label{eq3.8}&\geq
\tr (\bb E_{{\cl N}}(a)\bb E_{{\cl N}}(vv^*))
-\|b\|_2\|vv^*\|_2.
\end{align}
We now estimate the last term in (\ref{eq3.6}). By (\ref{eq3.5}),
\begin{align}
\|\bb
E_{\cl N}(ugv)\|^2_2 &= \|\bb E_{{\cl N}}(u)g\bb E_{{\cl
N}}(v)\|^2_2\nonumber\\
&= \tr(\bb E_{{\cl N}}(v)^*g^{-1}\bb E_{{\cl N}}(u)^*
\bb E_{{\cl N}}(u)g \bb E_{{\cl N}}(v))\nonumber\\
&\leq \tr(\bb E_{{\cl N}}(v)^*g^{-1}\bb E_{{\cl N}}(u^*u)g\bb E_{{\cl
N}}(v))\nonumber\\
&=\tr(\bb E_{{\cl N}}(a)\bb E_{{\cl N}}(v)\bb E_{{\cl N}}(v)^*)
+\tr(g^{-1}\bb E_{{\cl N}}(b)g\bb E_{{\cl
N}}(v)\bb E_{{\cl N}}(v)^*)\nonumber\\
\label{eq3.9}&\leq \tr(\bb E_{{\cl N}}(a)\bb E_{{\cl N}}(v)\bb E_{{\cl
N}}(v)^*) +\|\bb E_{{\cl N}}(v)\bb E_{{\cl N}}(v)^*\|_2\|b\|_2.
\end{align}
Using (\ref{eq3.8}) and (\ref{eq3.9}), (\ref{eq3.6}) becomes
\begin{align}
\|(I-\bb E_{\cl N})(u\bb E_{\cl N}(\cdot)v)\|^2_{\infty,2} &\ge 
\tr[\bb E_{\cl N}(a)(\bb E_{\cl N}(vv^*)
 -\bb E_{\cl N}(v)
\bb E_{\cl N}(v)^*)]\nonumber\\
\label{eq3.10} &\quad -\|b\|_2(\|vv^*\|_2+\|\bb E_{{\cl N}}(v)\bb E_{{\cl
N}}(v)^*\|_2).
\end{align}
Replacing $a$ and $b$ from (\ref{eq3.7}) gives (\ref{eq3.3}).

If $u$ is a unitary, (\ref{eq3.4}) follows immediately from (\ref{eq3.3}) by
replacing $u^*u$ with 1.
\end{proof}

\begin{thm}\label{thm3.3}
Let $G$ be an infinite subgroup of a countable discrete I.C.C. group $\Gamma$ 
with the property that $xGx^{-1}\cap G=\{e\}$ for all $x\in \Gamma \backslash
G$, let ${\cl M}=VN(\Gamma)$, and let ${\cl N}=VN(G)$.
\begin{itemize}
\item[(i)] If $u$ is a unitary in $\cl M$, then
\begin{equation}\label{eq3.11}
\|u-\bb E_{{\cl N}}(u)\|_2 \leq \|\bb E_{u{\cl N}u^*}-\bb E_{{\cl
N}}\|_{\infty,2};
\end{equation}
\item[(ii)] If $G$ is abelian, then $\cl N$ is a strongly singular masa
satisfying \begin{equation}\label{eq3.12}
\alpha(\cl N)=\delta(\cl N)=1.
\end{equation}\end{itemize}\end{thm}

\begin{proof}
\noindent (i)~~For $x \in {\cl M}$, $\|x\| \leq 1$, we have
\begin{align}
\|\bb E_{u\cl Nu^*} - \bb E_{\cl N}\|^2_{\infty,2} &\ge \|\bb E_{u\cl Nu^*}
(uxu^*) - \bb E_{\cl N}(uxu^*)\|^2_2\nonumber\\
&= \|u\bb E_{\cl N}(x) u^*-\bb E_{\cl N}(uxu^*)\|^2_2\nonumber\\
&\ge \|(I-\bb E_{\cl N}) [u\bb E_{\cl N}(x)u^* - \bb E_{\cl N}(uxu^*)]\|^2_2
\nonumber\\
\label{eq3.13}&= \|(I-\bb E_{\cl N})(u\bb E_{\cl N}(x)u^*)\|^2_2. 
\end{align}
Taking the supremum in (\ref{eq3.13}) over $x$ implies that
\begin{equation}\label{eq3.13a} \|\bb E_{u\cl Nu^*} - \bb E_{\cl N}\|^2_{\infty,2}
\geq \|(I-\bb E_{\cl N})(u\bb E_{\cl
N}(\cdot)u^*)\|^2_{\infty,2}. \end{equation}
Applying (\ref{eq3.4}) with $v=u^*$ gives
\begin{equation}\label{eq3.14}
\|\bb E_{u\cl Nu^*} - \bb E_{\cl N}\|^2_{\infty,2} \geq 1-\|\bb E_{\cl
N}(u)\|^2_2=\|u-\bb E_{\cl N}(u)\|^2_2,
\end{equation}
which proves (\ref{eq3.11}).

\noindent(ii)~~Assume now that $G$ is abelian. The estimate in  (\ref{eq3.11})
shows that $\cl N$ is a  strongly singular masa in $\cl M$ 
(although, {\it{a priori}}, it was not clear that $\cl N$ was
maximal). Thus, $\alpha(\cl N)=1$.

We now estimate $\delta(\cl A)$. Let $p$ and $q$ be orthogonal projections in
$\cl A$ which are equivalent in $VN(\Gamma)$. We may choose a nilpotent
partial isometry $v\in VN(\Gamma)$ such that $p = vv^*$ and $q=v^*v$. Then
\begin{equation}\label{eq3.16}
\bb E_{\cl N}(v) = \bb E_{\cl N}(pvq) = pq \bb E_{\cl N}(v) = 0.
\end{equation}
By (\ref{eq3.3}) and (\ref{eq3.16}),
\begin{equation}\label{eq3.17}
\|(I-\bb E_{\cl N}) v\bb E_{\cl N}(\cdot) v^*\|^2_{\infty,2} \ge \tr(q) =
\tr(p) = \|vv^*\|^2_2,
\end{equation}
since $v^*v \in {\cl N}'\cap {\cl M}={\cl N}$.
It follows that $\delta(\cl N)\ge 1$, and since $\delta(\cl N)\le 1$ is always
true, equality is immediate.

\end{proof}

We are now able to give some examples of strongly singular masas, which also
have the property that $\delta(\cl A)=1$.

\begin{cor}\label{cor3.4}
Let $\bb F_n$, $2\le n\le \infty$, denote the free group on $n$ generators,
let $a$ be one of these generators and let $\cl A$ be the masa generated by
$a$. Then $\cl A$ is strongly singular and $\delta(\cl A) = 1$.
\end{cor}

\begin{proof}
Any generator satisfies the hypotheses of Theorem \ref{thm3.3}.
\end{proof}

An element $a$ in a discrete group $\Gamma$ is prime, \cite{Gr,HV}, if the
equation $a=b^n$ has only two solutions in $\Gamma\colon \ b=a$ and $n=1$, or
$b=a^{-1}$ and $n=-1$. This says that $a$ is not a proper power of some other
group element.
The following lemma is surely well known, but we do not know a reference.

\begin{lem}\label{lem3.5} 
Let $a$ be a prime element of a group $\Gamma$. Then $xGp(a)x^{-1} \cap
Gp(a) = \{e\}$ for all $x\in \Gamma\backslash Gp(a)$ if and only if the
normalizer $N(Gp(a^p))$ of $Gp(a^p)$ is $Gp(a)$ for all $p\in \bb N$.
\end{lem}

\begin{proof}
One direction is clear. Conversely, suppose that the hypotheses on the
normalizers are fulfilled, but suppose that there is an $x\in \Gamma
\backslash Gp(a)$ such that, for some $p\in \bb N$ and $k\in \bb
Z\backslash\{0\}$, $xa^px^{-1} = a^k$. Then
\begin{equation}\label{eq3.18}
(x^{-1}ax)a^{pk}(x^{-1}a^{-1}x) = x^{-1}aa^{k^2}a^{-1}x = a^{pk},
\end{equation}
and so $x^{-1}ax \in N(Gp(a^{pk}))$. Hence $x^{-1}ax = a^r$ for some $r\in
\bb Z$, since one of $pk$, $-pk$ is in $\bb N$. Since $a$ is prime, so too is
$x^{-1}ax$, forcing $r=\pm 1$. Thus $x$ normalizes $Gp(a)$, a contradiction
which proves the result.
\end{proof}

Our next examples of strongly singular masas include those of Corollary
\ref{cor3.4}, and are based on a group theoretic result of Gromov, \cite{Gr}.

\begin{thm}\label{thm3.6}
Let $\cl A$ be the abelian von Neumann algebra generated by a prime element
$a$ in a non-elementary I.C.C. hyperbolic group $\Gamma$. Then $\cl A$ is a
strongly singular masa in $VN(\Gamma)$, and $\delta(\cl A)=1$.
\end{thm}

\begin{proof}
By \cite{Gr}, (see also Theorem~8.30 of \cite{GH}), a prime element $a$ in a
non--elementary hyperbolic group $\Gamma$ satisfies 
\begin{equation}\label{eq3.19}
N(Gp(a^p))=Gp(a)
\end{equation}
for all $p \in {\mathbb{N}}$. Lemma \ref{lem3.5} them shows that the
hypotheses of Theorem \ref{thm3.3} are satisfied, and the 
result follows. \end{proof}

\begin{cor}\label{cor3.7}
Let $n < m$, let ${\bb F}_m$ be the free group with generators
$\{g_i\}_{i=1}^m$, and regard ${\bb F}_n$ as a subgroup
generated by $\{g_i\}_{i=1}^n$. If ${\cl M}=VN({\bb F}_m)$ and
${\cl N}=VN({\bb F}_n)$, then
\begin{equation}\label{eq3.20}
\|u-\bb E_{{\cl N}}(u)\|_2 \leq \|\bb E_{u{\cl N}u^*}-\bb E_{{\cl
N}}\|_{\infty,2}
\end{equation}
for all unitaries $u \in \cl M$.
\end{cor}

\begin{proof}
The subgroup ${\bb F}_n$ of ${\bb F}_m$ satisfies the hypotheses of
Theorem \ref{thm3.3}.
\end{proof}

The following corollary provides an example of a strongly singular masa
$ \cl A$ in the hyperfinite type $II_1$ factor $\cl R$, and it also has Popa
invariant 1.

\begin{cor}\label{cor3.8}
In the hyperfinite  type $II_1$ factor $\cl R$, there exists a masa $\cl A$
satisfying
\begin{equation}\label{eq3.21}
\alpha(\cl A)=\delta(\cl A)=1.
\end{equation}
\end{cor}

\begin{proof}
Dixmier, \cite[Theorem 1]{Dix}, and Popa, \cite[Theorem 5.1]{Po3}, have
both given examples of countable amenable discrete I.C.C. groups containing
abelian subgroups which satisfy the hypotheses of Theorem \ref{thm3.3}, and
the result is then immediate.

For the reader's convenience, we briefly describe Dixmier's example. Let $K$
be an infinite field that is the countable union of finite subfields (the
algebraic closure of a finite field has this property). Let $\Gamma$ be the
group of affine transformations of the linear space of dimension 1 over $K$,
and let $G$ be the abelian subgroup of homotheties about 0. The calculations
of \cite[p.282]{Dix} show that  the hypotheses of Theorem \ref{thm3.3}
are satisfied.
\end{proof}

\begin{rem}\label{rem3.7}
Let $\cl A$ be a masa in a type $II_1$ factor $\cl M$, and let $\omega$ be a
free ultrafilter on $\mathbb N$. Then ${\cl A}^{\omega}$ is a masa in 
${\cl M}^{\omega}$ which is stongly singular when $\cl A$ also has this
property. The proof is similar to Popa's proof that 
$\delta({\cl A}^{\omega})=\delta(\cl A)$ for a masa $\cl A$ in  $\cl M$
\mbox{(\cite[Section 5.2]{Po2})}. There is also a version corresponding to a
sequence of strongly singular masas, again following 
 \mbox{\cite[Section 5.2]{Po2}}. $\hfill\square$ \end{rem}
\newpage

\section{Asymptotic homomorphism conditional expectations}\label{sec6}

\setcounter{equation}{0}

\indent

In this section we introduce the notion of an asymptotic homomorphism for
conditional expectations, and we show the certain abelian algebras arising
from group elements have this property. We then discuss some applications.

\begin{defn}\label{defn6.1}
Let $\cl A$ be an abelian von Neumann subalgebra of a type $II_1$ factor $\cl
M$. The conditional expectation $\bb E_{\cl A}$ is an asymptotic homomorphism
if there is a unitary $u\in\cl A$ such that
\begin{equation}\label{eq6.1}
\lim_{|k|\to\infty} \|\bb E_{\cl A}(xu^ky) - \bb E_{\cl A}(x) \bb E_{\cl A}(y)
u^k\|_2 = 0
\end{equation}
for all $x,y\in \cl M$.$\hfill\square$
\end{defn}

Observe that there is a closely related weak limit that converges for all
masas $\cl A$ in $\cl M$. Let $u$ be a unitary generating $\cl A$ and let
LIM be a Banach limit on $\bb Z$. Then, for $x,z\in \cl M$, we claim that
\begin{equation}\label{eq6.2}
\text{LIM}\left\langle n^{-1}  \sum^n_{j=1} u^{-j} xu^j,z\right\rangle =
\langle \bb E_{\cl A}(x),z\rangle.
\end{equation}
The left hand side of (\ref{eq6.2}) defines a bounded map $\phi\colon \ \cl
M\to \cl M$ by
\begin{equation}\label{eq6.3}
\text{LIM}\left\langle n^{-1} \sum^n_{j=1} u^{-j}xu^j,z\right\rangle = \langle
\phi(x),z\rangle,
\end{equation}
and the invariance of LIM shows that
\begin{align}
\langle u\phi(x),z\rangle &= \text{LIM}\left\langle n^{-1} \sum^n_1 u^{-(j-1)}
xu^{j-1}u,z\right\rangle\nonumber\\
\label{eq6.4}
&= \langle\phi(x)u,z\rangle.
\end{align}
Thus $\phi(x) \in \cl A'\cap \cl M=\cl A$ for all $x\in \cl M$, and since
$\phi$ is trace preserving, it is clear that $\phi=\bb E_{\cl A}$. Also, for
$x,y,z\in \cl M$, 
\begin{align}
&\text{LIM } n^{-1} \sum^n_{j=1} \langle u^{-j} \bb E_{\cl A}(xu^jy) - \bb
E_{\cl A}(x) \bb E_{\cl A}(y), z\rangle\nonumber\\
&\quad = \langle \bb E_{\cl A}(\bb E_{\cl A}(x)y) - \bb E_{\cl A}(x) \bb E_{\cl
A}(y), z\rangle\nonumber\\
\label{eq6.5}
&\quad = 0,
\end{align}
where we have used the fact that $\bb E_{\cl A}$ is both normal and an $\cl
A$-bimodule map.

\begin{thm}\label{thm6.2}
Let $\Gamma$ be a discrete group, let $\cl M = VN(\Gamma)$, and let $\cl A$ be
the abelian von~Neumann algebra generated by a fixed element $g\in\Gamma$. If
$g$ has the property that
$$\{k\in\bb Z\colon \ xg^k y\in Gp(g)\}$$
is finite for each pair $x,y\in\Gamma\backslash Gp(g)$, then $\bb E_{\cl A}$
is an asymptotic homomorphism.
\end{thm}

\begin{proof}
For each $k\in \bb Z$ define a bounded bilinear map $\phi_k\colon \ \cl
M\times \cl M\to \cl M$ by
\begin{equation}\label{eq6.6}
\phi_k(x,y) = \bb E_{\cl A}(xg^ky) - \bb E_{\cl A}(x) \bb E_{\cl A}(y) g^k
\end{equation}
for $x,y\in \cl M$. We consider first the case where $x$ and $y$ are group
elements in $\Gamma$. If either one is in $Gp(g)$ then the module properties
of $\bb E_{\cl A}$ imply that $\phi_k(x,y) = 0$ for all $k\in \bb Z$.

Now suppose that $x,y\in \Gamma \backslash Gp(g)$. By hypothesis, there
exists $K$ such that $xg^ky\notin Gp(g)$ for $|k|\ge K$, so both terms on the
right hand side of (\ref{eq6.6}) are 0, showing that $\phi_k(x,y) = 0$ for all
$|k|\ge K$. It then follows that, for $x,y\in \bb C\Gamma$, $\phi_k(x,y) = 0$
for $|k|$ sufficiently large.

The estimate
\begin{equation}\label{eq6.7}
\|\phi_k(x,y)\|_2 \le 2\|x\|_2 \|y\|
\end{equation}
for $x,y\in \cl M$, $k\in\bb Z$ is immediate from (\ref{eq6.6}), so if $x\in
\cl M$, $\{x_n\}^\infty_{n=1} \in \bb C\Gamma$, $y\in \bb C\Gamma$, and\break
$\lim\limits_{n\to \infty} \|x-x_n\|_2=0$, we obtain
\begin{align}
\|\phi_k(x,y)\|_2 &\le \|\phi_k(x-x_n,y)\|_2 + \|\phi_k(x_n,y)\|_2\nonumber\\
\label{eq6.8}
&\le 2\|x-x_n\|_2 \|y\| + \|\phi_k(x_n,y)\|_2.
\end{align}
Thus, for each $n\ge 1$,
\begin{equation}\label{eq6.9}
\mathop{\overline{\rm lim}}\limits_{|k|\to \infty} \|\phi_k(x,y)\|_2 \le
2\|x-x_n\|_2 \|y\|,
\end{equation}
since $\phi_k(x_n,y) = 0$ for $k$ sufficiently large. Let $n\to \infty$ in
(\ref{eq6.9}) to see that $\lim\limits_{|k|\to \infty} \|\phi_k(x,y)\|_2 = 0$
for $x\in \cl M$ and $y\in \bb C\Gamma$. Equation (\ref{eq6.6}) also gives the
estimate
\begin{equation}\label{eq6.10}
\|\phi_k(x,y)\|_2 \le 2\|x\| \ \|y\|_2.
\end{equation}
We then repeat the previous argument, this time in the second variable, to
obtain\hfil\break $\lim\limits_{|k|\to \infty} \|\phi_k(x,y)\|_2  =~0$ for all
$x,y\in \cl M$. This completes the proof.
\end{proof}

\begin{cor}\label{cor6.3}
Let $g$ be a prime element in a non-elementary I.C.C. hyperbolic group
$\Gamma$, and let $\cl A$ be the masa generated by $g$ in $VN(\Gamma)$. Then
$\bb E_{\cl A}$ is an asymptotic homomorphism.
\end{cor}

\begin{proof}
From the third section, $g$ satisfies the hypotheses of Theorem \ref{thm6.2},
and the result follows.
\end{proof}

\begin{rem}\label{rem6.4}
We remind the reader that Corollary \ref{cor6.3} applies, in particular, to
the generators of free groups.$\hfill\square$
\end{rem}

We now consider some consequences of asymptotic homomorphisms. We
will need the following inequality, which is close to Lemma \ref{lem3.2} under
different hypotheses.

\begin{pro}\label{pro6.5}
Let $\cl A$ be an abelian von~Neumann subalgebra of a type $II_1$ factor $\cl
M$. If $\bb E_{\cl A}$ is an asymptotic homomorphism, then
\begin{equation}\label{eq6.11}
\|(I-\bb E_{\cl A})(x\bb E_{\cl A}(\cdot)y)\|^2_{\infty,2} \ge \tr(\bb E_{\cl
A}(x^*x)\bb E_{\cl A}(yy^*) - \bb E_{\cl A}(x) \bb E_{\cl A}(x)^* \bb E_{\cl
A}(y) \bb E_{\cl A}(y)^*)
\end{equation}
for all $x,y\in \cl M$.
\end{pro}

\begin{proof}
In proving (\ref{eq6.11}), it clearly suffices to assume that $\|x\|$,
$\|y\|\le 1$. Now fix $\vp>0$. By the asymptotic homomorphism hypothesis we
may choose a unitary $u\in \cl A$ such that
\begin{align}
\label{eq6.12}
&\|\bb E_{\cl A}(xuy)  - \bb E_{\cl A}(x)  \bb E_{\cl A}(y)u\|_2 < \vp\\
\intertext{and}
\label{eq6.13}
&\|\bb E_{\cl A}(x^*xuyy^*) - \bb E_{\cl A}(x^*x) \bb E_{\cl A}(yy^*) u\|_2 <
\vp.
\end{align}
Then, using (\ref{eq6.12}) and (\ref{eq6.13}),
\begin{align}
\|(I-\bb E_{\cl A}) (x\bb E_{\cl A}(\cdot)y)\|^2_{\infty,2} &\ge \|(I-\bb
E_{\cl A})(xuy)\|^2_2\nonumber\\
&= \|xuy\|^2_2 - \|\bb E_{\cl
A}(xuy)\|^2_2\nonumber\\
&= \tr(x^*xuyy^*u^*) - \|\bb E_{\cl A}(xuy)\|^2_2\nonumber\\
&= \tr(\bb E_{\cl A}(x^*xuyy^*)u^*) - \|\bb E_{\cl A}(xuy)\|^2_2\nonumber\\
\label{eq6.14}
&\ge \tr(\bb E_{\cl A}(x^*x) \bb E_{\cl A}(yy^*))-\vp - \|\bb E_{\cl A}(x)
\bb E_{\cl A}(y)\|^2_2-2\vp.
\end{align}
Since $\bb E_{\cl A}(x) = \bb E_{\cl A}(y)$ commute,
\begin{equation}\label{eq6.15}
\|\bb E_{\cl A}(x) \bb E_{\cl A}(y)\|^2_2 = \tr(\bb E_{\cl A}(x) \bb E_{\cl
A}(x)^* \bb E_{\cl A}(y) \bb E_{\cl A}(y)^*),
\end{equation}
and so (\ref{eq6.11}) follows by substituting (\ref{eq6.15}) into
(\ref{eq6.14}) and letting $\vp\to 0$.
\end{proof}

\begin{rem}\label{rem6.6}
With only minor modifications to the proof, (\ref{eq6.11}) could be
strengthened \newline (under the same hypotheses) to 
\begin{align}
&\left\|(I-\bb E_{\cl A}) \left(\sum^n_{j=1} x_j \bb E_{\cl A}(\cdot)
y_j\right) \right\|^2_{\infty,2}\ge\nonumber\\
\label{eq6.16}
&\quad \tr\left(\sum^n_{i,j=1} [\bb E_{\cl A}(x^*_ix_j) \bb E_{\cl
A}(y_jy^*_i) - \bb E_{\cl A}(x_i)^* \bb E_{\cl A}(x_j) \bb E_{\cl A}(y_j) \bb
E_{\cl A}(y_i)^*]\right).
\end{align}
\end{rem}

Note that, in the particular case of a singly generated
subgroup, the following theorem gives the conclusions of Theorem \ref{thm3.3}
and  the asymptotic homomorphism condition appeared implicitly in the
proof of Lemma \ref{lem3.2}. 

\begin{thm}\label{thm6.7}
Let $\cl A$ be an abelian von Neumann subalgebra of a type $II_1$ factor
$\cl M$. If $\bb E_{\cl A}$ is an asymptotic homomorphism, then $\cl A$ is a
strongly singular masa with Popa invariant $\delta(\cl A)=1$.
\end{thm}

\begin{proof}
Let $u\in \cl M$ be an arbitrary unitary. To show strong singularity, we will
apply \newline Proposition~\ref{pro6.5} with $x=u$ and $y=u^*$. Then 
\begin{align}
\|\bb E_{u\cl A u^*} - \bb E_{\cl A}\|^2_{\infty,2} &= \|u\bb E_{\cl
A}(\cdot)u^* - \bb E_{\cl A}(u\cdot u^*)\|^2_{\infty,2}\nonumber\\
&\ge \|u\bb E_{\cl A}(\cdot)u^* - \bb E_{\cl A}(u\bb E_{\cl A}(\cdot)
u^*)\|^2_{\infty,2}\nonumber\\
&= \|(I-\bb E_{\cl A})(u\bb E_{\cl A}(\cdot) u^*)\|^2_{\infty,2}\nonumber\\
&\ge 1 - \|\bb E_{\cl A}(u) \bb E_{\cl A}(u^*)\|^2_2\nonumber\\
&\ge 1 - \|\bb E_{\cl A}(u)\|^2_2\nonumber\\
\label{eq6.17}
&= \|(I-\bb E_{\cl A})(u)\|^2_2.
\end{align}
Thus $\cl A$ is strongly singular. We now estimate $\delta(\cl A)$.

Let $v$ be  a nilpotent partial isometry such that $p=vv^*$ and $q=v^*v$ are
orthogonal projections in $\cl A$. From (\ref{eq3.16}), $\bb E_{\cl A}(v) = 0$,
so the choices of $x=v$ and $y=v^*$ in (\ref{eq6.11}) lead to
\begin{align}
\|(I-\bb E_{\cl A})(v\bb E_{\cl A}(\cdot)v^*)\|^2_2 &\ge \tr(\bb E_{\cl
A}(v^*v) \bb E_{\cl A}(v^*v))\nonumber\\
&= \tr(q)\nonumber\\
&= \tr(p)\nonumber\\
\label{eq6.18}
&= \|vv^*\|^2_2.
\end{align}
This proves that $\delta(\cl A)\ge 1$, and the reverse inequality always holds.
\end{proof}

\begin{rem}\label{rem6.8}
By inserting different powers of $u$ between elements of $\cl M$ in Definition \ref{defn6.1}
and letting the powers tend to $\infty$ {\it{successively}}, we can easily deduce 
a multivariable version as follows. If $\bb E_{\cl A}$ is an asymptotic homomorphism
and $x_1,\ldots,x_{n+1} \in \cl M$, then
\begin{equation}\label{eq6.19}
\lim_{k_1 \to \infty} \ldots \lim_{k_{n} \to \infty}
\|{\bb E_{\cl A}}(x_1u^{k_1}x_2\ldots u^{k_n}x_{n+1})-
{\bb E_{\cl A}}(x_1)\ldots {\bb E_{\cl A}}(x_{n+1})u^{k_1+\ldots +k_n}\|_2=0.
\end{equation}

If we insisted on using the same power in all $n$ places, then we would have
a type of freeness for $\bb E_{\cl A}$ and $\cl A$: for all $x_1,\ldots,x_{n+1} \in \cl M$,
\begin{equation}\label{eq6.20}
\lim_{|k| \to \infty}\|\bb E_{\cl A}(x_1u^kx_2\ldots u^kx_{n+1})-
\bb E_{\cl A}(x_1)\ldots \bb E_{\cl A}(x_{n+1})u^{nk}\|_2=0.
\end{equation}

Modifying the proof of Theorem \ref{thm6.2} shows that if $g$ is a generator of the free 
group ${\bb F}_n$ with associated masa $\cl A$, then $\bb E_{\cl A}$ has this freeness property using $u=g$ in
(\ref{eq6.20}). We have not investigated this idea. $\hfill\square$
\end{rem}

\newpage

\section{Strong singularity and the Popa invariant}\label{sec7}

\setcounter{equation}{0}

\indent

Recall that, for a masa $\cl A$ in a type $II_1$ factor $\cl M$, we defined
$\alpha(\cl A)$ to be the largest constant satisfying
\begin{equation}\label{eq7.1}
\|\bb E_{u\cl Au^*} - \bb E_{\cl A}\|_{\infty,2} \ge \alpha(\cl A)\|(I- \bb
E_{\cl A})(u)\|_2
\end{equation}
for all unitaries $u\in \cl M$. In this section we obtain an inequality
which links $\delta(\cl A)$ and $\alpha(\cl A)$. We will need the following
two lemmas.

\begin{lem}\label{lem7.1}
Let $\phi,\psi: \cl M \to \cl M$ be linear maps on a  type $II_1$ factor $\cl M$,
bounded in the $\|\cdot \|_{\infty,2}$--norm, and suppose that their ranges
are orthogonal in $L^2(\cl M, \tr)$. Then
\begin{equation}\label{eq7.1a}
\|\phi \pm \psi \|_{\infty,2} \leq \sqrt{\|\phi \|_{\infty,2}^2 +
\|\psi \|_{\infty,2}^2}\,.
\end{equation}\end{lem}

\begin{proof}
This is an immediate consequence of
\begin{equation}\label{eq7.1b}
\|h \pm k\| = \sqrt{\|h\|^2+\|k\|^2}
\end{equation}
for any pair of orthogonal vectors in a Hilbert space.
\end{proof}

\begin{lem}\label{thm5.3}
Let $\cl A$ and $\cl B$ be von Neumann subalgebras of a type $II_1$
factor $\cl M$. Then the following inequality holds:
$$\|\bb E_{\cl A'\cap \cl M}(I-\bb E_{\cl B'\cap\cl
M})\|_{\infty,2} \le 2\|(I-\bb E_{\cl A})\bb E_{\cl B}\|_{\infty,2}.$$
\end{lem}\
\begin{proof}
Let $x,y\in \cl M$ with $\|x\|, \|y\|_2\le 1$. Write $w = \bb E_{\cl
A'\cap \cl M}(y)$ and let $u\in \cl U(\cl B)$. Note that $w\in \cl A'$ and
$\|w\|_2 \le 1$. Then
\begin{align}
|\langle \bb E_{\cl A'\cap \cl M}(x-uxu^*),y\rangle| &= |\langle x-uxu^*, w\rangle|\nonumber\\
&= |\tr((x-uxu^*)w^*)|\nonumber\\
&= |\tr(xw^* - xu^*w^*u)|\nonumber\\
\label{eq5.15}
&\le |\tr(xw^* -xu^* \bb E_{\cl A}(u) w^*| + |\tr(xu^*w^*(u-\bb E_{\cl
A}(u)))|.
\end{align}
Here we have used the module properties of conditional expectations and that
$w^*$ and $\bb E_{\cl A}(u)$ commute.  The last expression in (\ref{eq5.15})
is no greater than
\begin{align}
&|\tr(xu^*(u-\bb E_{\cl A}(u)) w^*)| + |\tr(xu^*w^*(u-\bb E_{\cl
A}(u)))|\nonumber\\
&\quad \le 2\|(I-\bb E_{\cl A})(u)\|_2\nonumber\\
\label{eq5.16}
&\quad \le 2\|(I-\bb E_{\cl A})\bb E_{\cl B}\|_{\infty,2}
\end{align}
since $u = \bb E_{\cl B}(u)$. The estimates (\ref{eq5.15}) and (\ref{eq5.16})
combine to yield
\begin{equation}\label{eq5.17}
\|E_{\cl A'\cap \cl M}(x-uxu^*)\|_2\le 2\|(I-\bb E_{\cl A})\bb E_{\cl
B}\|_{\infty,2},
\end{equation}
letting $y$ vary over the unit ball of $L^2(\cl M, ~\tr)$. Since
\[\bb E_{\cl B'\cap \cl M}(x) \in
\overline{\text{conv}}^{\|\cdot\|_2}\{uwu^*\colon \ u\in \cl U(\cl B)\},\]
(see \cite{Ch,Po4}), the last inequality gives 
 \begin{equation}\label{eq5.18} \|E_{\cl A'\cap\cl M}(I-\bb E_{\cl
B'\cap \cl M})(x)\|_2 \le 2\|(I-\bb E_{\cl A})\bb E_{\cl B}\|_{\infty,2}.
\end{equation}
The result follows by letting $x$ vary over the unit ball of $\cl M$ in
(\ref{eq5.18}).
\end{proof}

\begin{thm}\label{thm7.1}
If $\cl A$ is a masa in a type $II_1$ factor $\cl M$, then
\begin{equation}\label{eq7.2}
\delta(\cl A) \ge \alpha(\cl A)/\sqrt{5}.
\end{equation}
In particular, $\delta(\cl A)\ge
1/\sqrt{5}$ for all strongly singular masas.
\end{thm}

\begin{proof}
Let $v$ be a nilpotent partial isometry such that $p=vv^*$ and $q=v^*v$ are
orthogonal projections in $\cl A$, and define a unitary $u\in \cl M$ by
\begin{equation}\label{eq7.3}
u=v+v^* + 1-p-q.
\end{equation}
Then
\begin{equation}\label{eq7.4}
\cl A = (1-p-q)\cl A + p\cl A + q\cl A,
\end{equation}
and
\begin{equation}\label{eq7.5}
u\cl Au^* = (1-p-q) \cl A + v\cl Av^* + v^*\cl Av.
\end{equation}
Thus
\begin{equation}\label{eq7.6}
\bb E_{\cl A} - \bb E_{u\cl Au^*} = \bb E_{p\cl A}+ \bb E_{q\cl A} - \bb
E_{v\cl Av^*} - \bb E_{v^*\!\cl Av}.
\end{equation}
By the orthogonality of $p$ and $q$, and the modularity of $\bb E_{\cl A}$,
\begin{equation}\label{eq7.7}
\|\bb E_{\cl A} - \bb E_{u\cl Au^*}\|^2_{\infty,2} = \|\bb E_{p\cl A} - \bb
E_{v\cl Av^*}\|^2_{\infty,2} + \|\bb E_{q\cl A} - \bb E_{v^*\!\cl
Av}\|^2_{\infty,2}.
\end{equation}
The two terms on the right hand side of (\ref{eq7.7}) are equal because the
map $v^*(\cdot)v$ implements an isometry from $p\cl Mp$ to $q\cl Mq$ in the
norms $\|\cdot\|$ and $\|\cdot\|_2$. Thus (\ref{eq7.7}) becomes
\begin{equation}\label{eq7.8}
\|\bb E_{\cl A}- \bb E_{u\cl Au^*}\|^2_{\infty,2}  = 2\|\bb E_{p\cl A} - \bb
E_{v\cl Av^*}\|^2_{\infty,2}.
\end{equation}
From (\ref{eq3.16}), $\bb E_{\cl A}(v) = 0$, and so the definition of $u$ gives
\begin{equation}\label{eq7.9}
\bb E_{\cl A}(u) = 1-p-q.
\end{equation}
Thus
\begin{equation}\label{eq7.10}
\|u-\bb E_{\cl A}(u)\|^2_2 = 1-\|\bb E_{\cl A}(u)\|^2_2
= \tr(p+q)
= 2 ~\tr(p)
= 2\|vv^*\|^2_2.
\end{equation}
Then (\ref{eq7.8}) and (\ref{eq7.10}) give
\begin{equation}\label{eq7.11}
\|\bb E_{p\cl A} - \bb E_{v\cl Av^*}\|_{\infty,2} \ge \alpha(\cl A)\|vv^*\|_2.
\end{equation}

In the von Neumann algebra $p\cl Mp$, consider the masas $p\cl A$ and $v\cl
Av^*$. We may apply Lemma~\ref{thm5.3}  to obtain
\begin{equation}\label{eq7.12}
\|\bb E_{p\cl A}(I-\bb E_{v\cl Av^*})\|_{\infty,2} \le 2\|(I-\bb E_{p\cl A})
\bb E_{v\cl Av^*}\|_{\infty,2}.
\end{equation}
Since $\bb E_{p\cl A}(I-\bb E_{v\cl Av^*})$ and 
$(I-\bb E_{p\cl A})\bb E_{v\cl Av^*}$ have orthogonal ranges,
\begin{align}
\|\bb E_{p\cl A} - \bb E_{v\cl Av^*}\|_{\infty,2} &= \|\bb E_{p\cl A}(I-\bb
E_{v\cl Av^*}) - (I-\bb E_{p\cl A})\bb E_{v\cl Av^*}\|_{\infty,2}\nonumber\\
\label{eq7.13}
&\le \sqrt{5}\|(I-\bb E_{p\cl A}) \bb E_{v\cl Av^*}\|_{\infty,2},
\end{align}
using (\ref{eq7.12}) and Lemma \ref{lem7.1}. Note that these expectations
are defined on $p\cl Mp$, but viewing them on $\cl M$ by first applying $\bb
E_{p\cl Mp}$ does not change the inequality (\ref{eq7.13}). We now combine
(\ref{eq7.11}) and (\ref{eq7.13}) to obtain
\begin{equation}\label{eq7.14}
\|(I-\bb E_{p\cl A}) \bb E_{v\cl Av^*}\|_{\infty,2} \ge  \alpha(\cl
A)\|vv^*\|_2/\sqrt{5},
\end{equation}
which shows that $\delta(\cl A) \ge \alpha(\cl A)/\sqrt{5}$.
\end{proof}

\begin{rem}\label{rem7.2}
Theorem \ref{thm7.1} raises several obvious questions. Can the factor of 
$\sqrt{5}$ be
removed from the inequality (\ref{eq7.2}) in this theorem? Is there an
inequality in the opposite direction? Is it possible that $\delta(\cl A) =
\alpha(\cl A)$ in general? $\hfill\square$ \end{rem}

\end{document}